\documentclass[a4paper,10pt]{article}
\usepackage{amsmath,amsthm,amssymb,amscd}

\newtheorem{theorem}{Theorem}[section]
\newtheorem{lemma}[theorem]{Lemma}
\newtheorem{proposition}[theorem]{Proposition}
\newtheorem{corollary}[theorem]{Corollary}

\theoremstyle{definition}

\theoremstyle{remark}

\numberwithin{equation}{section}

\begin{document}

\title{\bf Integral Transform
and\\ Segal-Bargmann Representation\\
associated to q-Charlier Polynomials\footnote
{Posted on \tt http://arXiv.org/abs/math.CA/0104260}}
\author{Nobuhiro ASAI\thanks{Research supported by a Postdoctral 
Fellowship of the International Institute for Advanced Studies}\\
International Institute for Advanced Studies\\
Kizu, Kyoto, 619-0225, Japan.\\
{\rm asai@iias.or.jp, nobuhiro.asai@nifty.com}}
\date{(\small Submitted April 24, 2001.)}
\maketitle

\begin{abstract}
{Let $\mu_p^{(q)}$ be the $q$-deformed 
Poisson measure in the sense of Saitoh-Yoshida \cite{sy2}
and $\nu_p$ be the measure given by Equation \eqref{eq:nu-q}.
In this short paper, 
we introduce the $q$-deformed analogue of 
the Segal-Bargmann transform associated with $\mu_p^{(q)}$.
We prove that our Segal-Bargmann transform 
is a unitary map of $L^2(\mu_p^{(q)})$ 
onto the $q$-deformed Hardy space ${\cal H}^2(\nu_q)$.
Moreover, we give the 
Segal-Bargmann representation of 
the multiplication operator by $x$ in $L^2(\mu_p^{(q)})$,
which is a linear combination of the $q$-creation, 
$q$-annihilation, $q$-number, and scalar operators.}
\end{abstract}

\section{Introduction}

The classical Segal-Bargmann transform 
in Gaussian analysis
yields a unitary map of $L^2$ space 
of the Gaussian measure on ${\Bbb R}$ 
onto the space of $L^2$ holomorphic functions of the
Gaussian measure on ${\Bbb C}$, 
see papers \cite{barg1,barg2,dwyer,gm,kk,segal1,segal2}.

Recently, Accardi-Bo\.{z}ejko \cite{ab}
showed the existence of a unitary operator between a one-mode
interacting Fock space and $L^2$ space of a probability measure
on ${\Bbb R}$ by making use of the basic properties of 
classical orthogonal polynomials 
and associated recurrence formulas \cite{chihara,szego}.
Inspired by this work, the author \cite{asai01} has recently 
extended the Segal-Bargmann transform to non-Gaussian 
cases.  The crucial point is to introduce a coherent state 
vector as a kernel function 
in such a way that a transformed function,
which is a holomorphic function on a certain domain in general, 
becomes a power series expression.  
Along this line, 
Asai-Kubo-Kuo \cite{akk01} 
have considered the case of 
the Poisson measure compared with 
the case of the Gaussian measure.
However, the case of 
$L^2$ space of Wigner's semi-circle distributions
in free probability theory \cite{vdn} is beyond their scope. 

On the other hand, 
Leeuwen-Maassen \cite{lm95} considered a transform 
associated with $q$ deformation of 
the Gaussian measure \cite{bs1,bs2,bks}
and showed that for a given real number $q\in [0,1)$
it is a unitary map of $L^2$ space of 
$q$-defomed Gaussian measure onto the $q$-deformed
Hardy space ${\cal H}^2(\nu_q)$
where $\nu_q$ is given in \eqref{eq:nu-q}.
Biane \cite{biane97} examined the case of 
$q=0$ (Free case).  Roughly speaking,  
their methods do not give the relationship between 
Szeg\"o-Jacobi parameters and kernel functions for 
their transforms.  
As observed in Section 
\ref{sec:transform} and Appendix \ref{sec:appendixA},
our approach clarifies the relationship between them.

In this paper, 
we shall consider the $q$-deformed version of 
the Segal-Bargmann transform $S_{\mu_p^{(q)}}$ 
associated with $q$-deformed 
Poisson measure, denoted by $\mu_p^{(q)}$,
in the sense of Saitoh-Yoshida \cite{sy2}.
As a main result, we shall provide Proposition \ref{prop:main}, 
which claims that $S_{\mu_p^{(q)}}$
is a unitary map of $L^2(\mu_p^{(q)})$ 
onto ${\cal H}^2(\nu_q)$.
Moreover, in Theorem \ref{thm:Q-p}
we shall give the 
representation in ${\cal H}^2(\nu_q)$ of 
the multiplication operator by $x$ in $L^2(\mu_p^{(q)})$,
which is a linear combination of the $q$-creation, 
$q$-annihilation, $q$-number, and scalar operators.
We remark that our representation is compatible with 
that on the $q$-Fock space by Saitoh-Yoshida \cite{sy1} 
and can be viewed as the $q$-analogue of 
the Hudson-Parthasarathy \cite{hp}
decomposition of the usual Poisson random 
variable on the standard Boson Fock space ($q=1$).
Ito-Kubo \cite{ik} also studied a similar decomposition
in details from the point of white noise calculus \cite{kt,kuo96}
(For more recent formulation, see papers \cite{akk6,akk7,cks}).

The present article serves a good example 
to papers by Accardi-Bo\.{z}ejko \cite{ab}
and Asai \cite{asai01}.
The present paper is organized as follows.
In Section \ref{sec:q}, we recall the recurrence formula 
for $q$-Charlier polynomials.  
In Section \ref{sec:transform},
we introduce a $q$-deformed coherent state vector and a
Segal-Bargmann transform associated 
to a $q$-deformed Poisson measure.
In addition, we quickly define the Hardy space as 
the Segal-Bargmann representation space.
In Section \ref{sec:main}, our main results are given.
In Appendix \ref{sec:appendixA}, 
we give some remarks on known 
results \cite{barg1,biane97,lm95,segal1,segal2}
related to $q$-Hermite polynomials.

\medskip
\noindent
{\bf Notation.}\\
Let us recall standard notation from $q$-analysis \cite{aar}.
We put for $n\in {\Bbb N_0}$, 
\begin{equation*}
	[n]_q:=1+q+\cdots +q^{n-1} \ \text{with} \ [0]_q=0. 
\end{equation*}
Then $q$-factorial is naturally defined as 
\begin{equation*}
	[n]_q!:=[1]_q\cdots [n]_q \ \text{with} \ [0]_q!=0. 
\end{equation*}
The $q$-exponential is given by 
\begin{equation*}
	\exp_q(x):=\sum_{n=0}^{\infty}{x^n\over [n]_q!}.
\end{equation*}
whose radious of convergence is $1/(1-q)$.
In addition, another symbol used is 
the q-analogue of the Pochhammer symbol,
\begin{equation*}
	(a;q)_n=\prod_{j=0}^{n-1}(1-aq^j) \ \text{and} \
	(a;q)_{\infty}=\prod_{j=0}^{\infty}(1-aq^j)
\end{equation*}
with the convention $(a;q)_0=1$.

\section{q-deformed Charlier Polynomials}\label{sec:q}
From now on, we always assume that $q\in[0,1)$ is fixed.
Recently, Saitoh-Yoshida \cite{sy2} calculated the 
explicit form of the $q$-deformed Poisson measure with
a parameter $\beta>0$ for $q\in [0,1)$.
We denote it by $\mu_{p}^{(q)}$.  
The orthogonal polynomials
associated to $\mu_p^{(q)}$ are
the $q$-Charlier polynomials $\{C_n^{(q)}(x)\}$ 
with the Szeg\"{o}-Jacobi parameters 
$\alpha_n=[n]_q+\beta, \omega_n=\beta [n]_q$. 
See papers \cite{sy1,sy2}.  We also refer the 
book \cite{chihara} for 
the standard Charlier polynomials case ($q=1$).    
Let $\lambda=\{\beta^n[n]_q!\}_{n=0}^{\infty}$.
The following relations hold for each $n\geq 1$:
\begin{eqnarray}\label{eq:Crecurrence}
	(x-[n]_q-\beta)C_n^{(q)}(x)=C_{n+1}^{(q)}(x)
	+\beta [n]_qC_{n-1}^{(q)}(x)
\end{eqnarray} 
where $C_{0}^{(q)}=1$, $C_1^{(q)}=x-\beta$.
Then, for any $L^2$-convergent decompositions
$f(x)=\sum a_nC_n^{(q)}(x)$ and 
$g(x)=\sum b_nC_n^{(q)}(x)$, the inner product 
$\langle\cdot,\cdot\rangle_{L^2(\mu_p^{(q)})}$ is given  
by the form
\begin{equation}
	\langle f, g\rangle_{L^2(\mu_p^{(q)})}
	=\sum_{n}\beta^n[n]_q!\overline{a}_nb_n.
\end{equation}
Moreover, let $\|f\|^2_{L^2(\mu_p^{(q)})}
=\sum_{n=0}^{\infty}\beta^n[n]_q!|a_n|^2$.

\section{Segal-Bargmann Transform and Hardy Space}
\label{sec:transform}
Let us define 
the {\it q-deformed coherent state vector} 
for $\{C_n^{(q)}(x)\}$ by 
\begin{equation}
	E^{(q)}_{\lambda,p}(x,z)=\sum_{n=0}^{\infty}
	{C_n^{(q)}(x)\over \beta^n[n]_q!}z^n,
	\quad z\in 
	\Omega^{\beta}_{q}:=
	\Bigl\{z\in {\Bbb C}:\ |z|^2<{\beta\over 1-q}\Bigr\}.
\end{equation}
It is easy to see that $E_p^{(q)}(x,z)\in 
L^2(\mu_p^{(q)})$ due to 
\begin{equation}
\|E^{(q)}_{\lambda,p}(x,z)\|^2_{L^2(\mu_p^{(q)})}
=\exp_q(\beta^{-1}|z|^2)<\infty.
\end{equation}
Moreover, it can be shown that 
$\{E^{(q)}_{\lambda,p}(x,y)\}$
is linearly independent and total in $L^2(\mu_p)$.

Now we are in a position to introduce our key tool 
in this paper.
Let us consider the $q$-analogue of the Segal-Bargmann 
transform $S_{\mu_p^{(q)}}$ associated to $\{C_n^{(q)}(x)\}$ 
given by 
\begin{equation}
	(S_{\mu_p^{(q)}}f)(z)
	=\langle E_{\lambda,p}^{(q)}(x,\overline{z}), 
	f(x)\rangle_{L^2(\mu_p^{(q)})} \ 
	\text{for any} \ f\in L^2(\mu_p^{(q)}).
\end{equation}

\begin{lemma}\label{lem:1}
Let $f\in L^2(\mu_p^{(q)})$.  
Then $(S_{\mu_p^{(q)}}f)(z)$ 
converges absolutely for all $z\in\Omega^{\beta}_q$.
\end{lemma}
\begin{proof}
For $f(x)=\sum_{n=0}^{\infty}a_nC_n^{(q)}(x)$, it 
is quite easy to see 
\begin{equation}
	(S_{\mu_p^{(q)}}f)(z)
	=\sum_{n=0}^{\infty}a_nz^n.
\end{equation}
By the Schwartz inequality, we get the inequality
\begin{equation*}
	\sum_{n=0}^{\infty}|a_nz^n|
	\leq \|f\|_{L^2(\mu_p^{(q)})}
	\exp_q\Bigl({|z|^2\over \beta}\Bigr).
\end{equation*}
This shows that the $(S_{\mu_p^{(q)}}f)(z)$ 
converges absolutely for all $z\in\Omega_q^{\beta}$.
\end{proof}

The completion of the space of holomorphic functions 
$F, G$ on $\Omega^{\beta}_{q}$ with respect to the inner product,
\begin{equation}
	(F, G)_{{\cal H}^2_q}
	=\int \overline{F(z)}
	G(z)\nu_q(dz),
\end{equation} 
is nothing but 
the $q$-deformed Hardy space ${\cal H}^2({\nu}_q)$.
Here ${\nu}_q$ means that 
\begin{equation}\label{eq:nu-q}
	\nu_q(dz)
	=(q;q)_{\infty}\sum_{j=0}^{\infty}{q^j\over (q;q)_j}
	\lambda^{\beta}_{r_j}(dz),\quad
	q\in [0,1), \ r_j=q^{j\over 2}{\sqrt{\beta\over 1-q}}
\end{equation}
where $\lambda^{\beta}_{r_j}(dz)$ is the Lebesgue measure 
on the circle of radious $r_j$.

\begin{lemma}\label{lem:2}
$\{z^n\}_{n=0}^{\infty}$ forms an orthogonal basis of 
${\cal H}^2(\nu_q)$.
\end{lemma}
\begin{proof}
We adopt the same idea as 
in the proof \cite{lm95}.
\begin{align*}
	(z^n, z^m)_{{\cal H}^2_q}
	& = \int_{\Omega_q^{\beta}}\overline{z}^nz^m\nu_q(dz)\\
	& =(q;q)_{\infty}\sum_{j=0}^{\infty}{q^j\over (q;q)_j}
	\int_{\Omega_q^{\beta}}
	\overline{z}^nz^m\lambda_{r_j}^{\beta}(dz)\\
	& ={(q;q)_{\infty}\over 2\pi}
	\sum_{j=0}^{\infty}{q^j\over (q;q)_j}r_j^{n+m}
	\int_{0}^{2\pi}e^{i(m-n)\theta}d\theta\\
	& =\delta_{n,m}{\beta^n(q;q)_{\infty}\over (1-q)^n}
	\sum_{j=0}^{\infty}{q^{(n+1)j}\over (q;q)_j}\\
	& =\delta_{n,m}\beta^n[n]_q!.
\end{align*}
Note that we have used the $q$-Gamma function \cite{aar,lm95},
\begin{equation*}
	\Gamma_q(n+1):={(q;q)_{\infty}\over (1-q)^n}
	\sum_{j=0}^{\infty}{q^{(n+1)j}\over (q;q)_j}
	=[n]_q! \ \text{for any} \ n\in {\Bbb N}.
\end{equation*}
\end{proof}

\medskip
\noindent
{\it Remark.}
It can be shown by Proposition 4.4 in the recent paper \cite{akk01}
that $\nu_q$ is a unique measure satisfying 
$(z^n, z^m)_{{\cal H}^2_q}=\delta_{n,m}\beta^n[n]_q!$.

\medskip
Hence, for any $F=\sum_{n=0}^{\infty}a_nz^n$, 
$G=\sum_{n=0}^{\infty}
b_nz^n\in {\cal H}^2({\nu}_q)$,
the inner product $(\cdot,\cdot)_{{\cal H}^2_q}$
is written as 
\begin{equation}
	(F,G)_{{\cal H}^2_q}
	=\sum_{n=0}^{\infty}\beta^n[n]_q!\overline{a}_nb_n
\end{equation}
and the corresponding norm of $F$ is 
\begin{equation}
	\|F\|^2_{{\cal H}^2_q}=\sum_{n=0}^{\infty}
	\beta^n[n]_q!|a_n|^2.
\end{equation}

\section{Main Results}\label{sec:main}
\begin{proposition}\label{prop:main}
$S_{\mu_p^{(q)}}$ is a unitary map of $L^2(\mu_p^{(q)})$
onto ${\cal H}^2({\nu}_q)$.
\end{proposition}
\begin{proof}
As we have seen in Lemma \ref{lem:1},
\begin{equation}\label{eq:z^n}
	(S_{\mu_p^{(q)}}C_n^{(q)})(z)=z^n.
\end{equation}
In addition, we derive by Lemma \ref{lem:2}, 
\begin{equation*}
	\|C_n^{(q)}\|^2_{L^2(\mu_p^{(q)})}
	=\|z^n\|^2_{{\cal H}^2_q}=\beta^n.
\end{equation*}
Therefore, we finish the proof.
\end{proof}

Let us define operators $Z_q$ and $D_q$ in 
${{\cal H}^2({\nu}_q)}$ satisfying
\begin{equation}\label{eq:q-creation}
	Z_qF(z)=zF(z).
\end{equation}
and 
\begin{equation}\label{eq:q-annihilation}
	D_{q,\beta}F(z):=0 \ (z=0), 
	\quad D_{q,\beta}F(z):={\beta(F(z)-F(qz))\over z(1-q)}.
	\ (z\ne 0)
\end{equation}
Operators $Z_q$ and $D_{q,\beta}$ play the roles of
the {\it q-creation operator} and {\it q-annihilation operator} 
respectively
and satisfy the $q$-deformed commutation relation
$D_{q,\beta}Z_q-qZ_qD_{q,\beta}=I$.
The {\it q-number operator} acting on 
${{\cal H}^2({\nu}_q)}$ is defined by 
\begin{equation}\label{eq:q-number}
	\widetilde{N}_qF(z)=
	[n]_qF(z), \ n\geq 0.
\end{equation}
In addition, the operator $\widetilde{\alpha}_{N_q}$ 
acting on 
${{\cal H}^2({\nu}_q)}$ is defined by 
\begin{equation}\label{eq:q-alpha}
	\widetilde{\alpha}_{N_q}F(z)=
	([n]_q+\beta)F(z), \ n\geq 0.
\end{equation}
Remark that $\widetilde{\alpha}_{N_q}F(z)=0$ for $q$-Gaussian case,
see Appendix \ref{sec:appendixA}.
By the direct calculation, we have 
\begin{lemma}\label{lem:}
\begin{itemize}
\item[$(1)$]
	$S_{\mu_p^{(q)}}1=1$
\item[$(2)$]
	$D_{q,\beta}z^n
	=\beta [n]_qz^{n-1}$
\item[$(3)$]
	$Z_qz^n
	=z^{n+1}$
\end{itemize}
\end{lemma}

The transformation of the multiplication 
operator $Q_p^{(q)}$ by $x$ in $L^2(\mu_p^{(q)})$
satisfies the following relation.
\begin{theorem}\label{thm:Q-p}
	$S_{\mu_p^{(q)}}Q_p^{(q)}
	=(D_{q,\beta}+Z_q+\widetilde{\alpha}_{N_q})S_{\mu_p^{(q)}}$
\end{theorem}
\begin{proof}
By the recurrence formula \eqref{eq:Crecurrence},
Equation \eqref{eq:z^n} and Lemma \ref{lem:}, we derive
\begin{align*}
	(S_{\mu_p^{(q)}}Q_p^{(q)}C_n^{(q)})(z)
	& =\bigl\langle E^{(q)}_{\lambda,p}(x,\overline{z}),
	C_{n+1}^{(q)}+\beta [n]_qC_{n-1}^{(q)}
	+([n]_q+\beta)C_n^{(q)}
	\bigr\rangle_{L^2(\mu^{(q)}_p)}\\
	& =(D_{q,\beta}+Z_q+\widetilde{\alpha}_{N_q})
	(S_{\mu_p^{(q)}}C_n^{(q)})(z).
\end{align*}
\end{proof}
Moreover we obtain
\begin{corollary}\label{cor:decomposition}
The operators $Z_q, D_q, \widetilde{N}_q, \widetilde{\alpha}_{N_q}$
have the following properties:\\
(1) $\widetilde{\alpha}_{N_q}
={1\over\beta}Z_qD_{q,\beta}+\beta I
=\widetilde{N}_q+\beta I$\\
(2) $Z_q+D_{q,\beta}+\widetilde{\alpha}_{N_q}=
	\bigl({1\over\sqrt{\beta}}Z_q+\sqrt{\beta}I\bigr)
	\bigl({1\over\sqrt{\beta}}D_{q,\beta}+\sqrt{\beta}I\bigr)$.
\end{corollary}
\begin{proof}
The proof is done by Equations \eqref{eq:q-creation},
\eqref{eq:q-annihilation}, \eqref{eq:q-number},
\eqref{eq:q-alpha}, and Lemma \ref{lem:}.
\end{proof}

Therefore, by Theorem \ref{thm:Q-p} 
and Corollary \ref{cor:decomposition},  
the multiplication operator $Q_p^{(q)}$ by $x$ 
in $L^2(\mu_p^{(q)})$ is represented as a linear
combination of the $q$-creation, 
$q$-annihilation, $q$-number, and scalar operators
in the Hardy space ${\cal H}^2(\nu_q)$.
In the $q$-deformed Gaussian case,
due to Equation \eqref{eq:Q-g} in Appendix \ref{sec:appendixA}, 
the multiplication operator $Q_g^{(q)}$ in $L^2(\mu_g^{(q)})$ 
is represented as a linear combination of $q$-creation and 
$q$-annihilation operators in the ``same" Hardy space 
${\cal H}^2(\nu_q)$.

\medskip
\noindent
{\it Remark.}
(1) The analogous results in this paper to $q=1$
have been considered by 
Asai, et al. \cite{akk01}.  
(2) In general, 
if variances of two given measures $\mu_1$ and 
$\mu_2$ are the same, 
then the actions of creation and annihilation
operators are the same.  In addition, the transformed
function spaces by $S_{\mu_1}$ and $S_{\mu_2}$
are also the same.
However,  if $\mu_1$ is non-symmetric
and $\mu_2$ is symmetric,
then the representation of multiplication in 
$L^2(\mu_1)$ is different from 
that in $L^2(\mu_2)$. 

\appendix
\section{Appendix. On q-Hermite Polynomials}\label{sec:appendixA}
\label{q-polynomial}
First of all, we refer to the papers \cite{bks,lm95} 
and references cited therein for the detailed description of 
$q$-Hermite polynomials.  For $q=1$, 
see books \cite{chihara,kuo96,szego}.

Let $\mu_{g}^{(q)}$ be
the $q$-deformed Gaussian measure with 
mean zero and variance $\beta>0$.  
It is well-known that an associated orthogonal
polynomial to $\mu_{g}^{(q)}$ is the $q$-Hermite polynomial
$\{H_n^{(q)}(x)\}$ with the Szeg\"{o}-Jacobi parameters $\alpha_n=0, 
\omega_n=\beta [n]_q$.  
In this case 
the following relations hold for each $n\geq 1$:
\begin{eqnarray}\label{eq:Hrecurrence}
	xH_n^{(q)}(x)=H_{n+1}^{(q)}(x)+\beta [n]_qH_{n-1}^{(q)}(x)
\end{eqnarray} 
where $H_{0}^{(q)}=1$ and $H_1^{(q)}=x$.
Then for any $f(x)=\sum a_nH_n^{(q)}(x)$ and 
$g(x)=\sum b_nH_n^{(q)}(x)$ in $L^2(\mu_g^{(q)})$, 
the inner product 
$\langle\cdot,\cdot\rangle_{L^2(\mu_g^{(q)})}$ is given  
by the form
\begin{equation*}
	\langle f, g\rangle_{L^2(\mu_g^{(q)})}
	=\sum_{n}\beta^n[n]_q!\overline{a}_nb_n.
\end{equation*}
Moreover, let $\|f\|^2_{L^2(\mu_g^{(q)})}
=\sum_{n=0}^{\infty}\beta^n[n]_q!|a_n|^2$.

A {\it q-deformed coherent state vector} 
$E^{(q)}_{\lambda,g}(x,z)$ 
for $\{H^{(q)}_n(x)\}$ is defined by 
\begin{equation}\label{eq:Eg-vector}
	E^{(q)}_{\lambda,g}(x,z)=\sum_{n=0}^{\infty}
	\frac{H^{(q)}_{n}(x)}{\beta^nn!}z^n,\quad
	z\in\Omega_q^{\beta}.
\end{equation}
It can be shown that the set
$\{E^{(q)}_{\lambda,g}(x,z): z\in \Omega_q^{\beta}\}$
is linearly independent and total in $L^2(\mu^{(q)}_g)$.
The $q$-analogue of the Segal-Bargmann 
transform $S_{\mu_g^{(q)}}$ associated to $\{H_n^{(q)}(x)\}$ 
is given by 
\begin{equation}
	(S_{\mu_g^{(q)}}f)(z)
	=\langle E_{\lambda,g}^{(q)}(x,\overline{z}), 
	f(x)\rangle_{L^2(\mu_g^{(q)})} \ 
	\text{for any} \ f\in L^2(\mu_g^{(q)}).
\end{equation}

With this transform, we can reproduce the same results 
as for $q\in[0,1)$ Theorem III.4 
in Leeuwen-Maassen \cite{lm95} and 
for $q=0$ Proposition 1
in Biane \cite{biane97}.
That is, 
$S_{\mu_g^{(q)}}$-transform yields a unitary isomorphism
between $L^2(\mu_g^{(q)})$ and ${\cal H}^2(\nu_q)$
and 
\begin{equation}\label{eq:Q-g}
	S_{\mu_g^{(q)}}Q_g^{(q)}
	=(D_{q,\beta}+Z_q)S_{\mu_g^{(q)}}
\end{equation}
where $Q_g^{(q)}$ is the multiplication operator by $x$
in $L^2(\mu_g^{(q)})$.  

\section*{Acknowledgments}
The author expresses sincere thanks to Professors
T. Hida and K. Sait\^o for their kind arrangement of the
conference.  He is grateful for research support by
International Institute for Advanced Studies, Kyoto, Japan.
He also thanks Professor Bo\.{z}ejko for his comments.

\end{document}